\documentclass{article}
\usepackage{amsmath}
\usepackage{amssymb}
\title{Some remarks on $A_1^{(1)}$ soliton cellular automata}
\author{Susumu Ariki}
\date{}
\begin{document}
\maketitle

\newtheorem{thm}{Theorem}
\newtheorem{cor}[thm]{Corollary}
\newtheorem{prop}[thm]{Proposition}
\newtheorem{defn}[thm]{Definition}
\newtheorem{examp}[thm]{Example}
\newtheorem{lem}[thm]{Lemma}

\newcommand{\BM}{\left[\begin{array}}
\newcommand{\EM}{\end{array}\right]}
\newcommand{\QED}{$\blacksquare$}
\newcommand{\Prod}{\Pi}
\newcommand{\Beta}{{\bf\beta}}

\begin{abstract}
We describe the $A_1^{(1)}$ soliton cellular 
automata as an evolution of a poset. 
This allows us to explain the conservation laws 
for the $A_1^{(1)}$ soliton cellular 
automata, one given by Torii, Takahashi and Satsuma, and 
the other given by 
Fukuda, Okado and Yamada, in terms of the 
stack permutations of states in a very natural manner. 
As a biproduct, we can prove a conjectured formula 
relating these laws. 
\end{abstract}

\section{Introduction}

Several years ago, Torii, Takahashi and Satsuma \cite{TTS} 
proved a conservation law for their box-ball system 
(soliton cellular automata) using 
the Robinson-Schensted-Knuth correspondence: 
we associate a permutation to each state $p$, 
which we call 
the stack permutation of the state, then the shape of the 
$P$-symbols of these stack permutations is conserved. 
We denote this partition by $\lambda(p)$. 

Recently, it was observed that there exists a crystal 
structure behind this system, and 
the identification of this box-ball system 
with a box-ball system arising from 
$A_1^{(1)}$-crystal was made in \cite{HHIKTT}. In this crystal 
picture, the time evolution 
is described by combinatorial 
row-to-row transfer matrices, and the energy 
functions $E_l(p)$ $(l\in{\mathbb N})$ of 
this system naturally gives us another 
conservation law \cite{FOY}. Further, 
it was conjectured 
how these laws were related. It is given by a simple formula: 
\[
E_l(p)-E_{l\!-\!1}(p)=\lambda_l(p)
\]
where $\lambda_l(p)$ is the length of the $l$ th row of 
the partition $\lambda(p)$. 

In terms of the lengths of solitons $N_1,N_2,\dots$, 
$\lambda(p)$ is the partition which has $N_k$ columns of 
length $k$ $(k\in{\mathbb N})$, and 
$E_l(p)=\sum_{k\in{\mathbb N}} {\rm min}(l,k)N_k$. 

For example, if the state is an asymptotic soliton state, 
\[
\cdots 01^{k_1}0^{l_1}1^{k_2}0^{l_2}\cdots\qquad(l_1,l_2,\cdots>\!>0)
\]
it is straightforward to verify it. 

The purpose of this 
short note is to prove the formula. It is done by 
supplying conceptual explanation about the appearence of 
stack permutations and their $P$-symbols. 

Main idea is to interprete the box-ball system as a 
discrete dynamical system of a path on 
${\mathbb Z}\times{\mathbb N}$, from which 
we naturally read off the evolution of the permutation poset 
of the stack permutation of the state of the 
original box-ball system. This gives us a natural 
explanation why the Torii-Takahashi-Satsuma law holds. 

We then turn to the crystal picture, and describe the 
sites which contribute to the energy function by using 
stack permutations. This explains why these energy functions 
are related to stack permutations. 

Our conclusion is that the depth of stacks explains 
both conservation laws, which proves the relation of 
these laws. 

The author hopes that this explanation would be valid after 
modifications in the case of $A_r^{(1)}$ soliton cellular 
automata. In this case, Nagai's conserved quantities 
remain mysterious from combinatorial point of view. 

\section{Fomin-Greene theory on posets}

We start with the Fomin-Greene theory of posets. 
Good references are \cite{BF} and \cite{F}. Let $({\cal P},\le)$ be 
a poset. A {\bf chain} is a totally ordered subset of 
${\cal P}$. An {\bf antichain} is a subset of ${\cal P}$ 
on which no two elements are comparable. 

\begin{defn}
Let $({\cal P},\le)$ be a poset. 
We define $I_k({\cal P}),\,D_k({\cal P})$ for $k\in{\mathbb N}$ 
as follows. 
\[
\begin{array}{c}
I_k({\cal P}):={\rm max}\{\,
|C_1\sqcup\cdots\sqcup C_k|\,|\,\mbox{
$C_i$:(possibly empty) chain}\,\}\\
D_k({\cal P}):={\rm max}\{\,
|A_1\sqcup\cdots\sqcup A_k|\,|\,\mbox{
$A_i$:(possibly empty) antichain}\,\}
\end{array}
\]
We also define $\lambda_k({\cal P}),\,\lambda'_k({\cal P})$ 
for $k\in{\mathbb N}$ by their differences: 
\[
\lambda_k({\cal P}):=I_k({\cal P})-I_{k\!-\!1}({\cal P}),\quad
\lambda'_k({\cal P}):=D_k({\cal P})-D_{k\!-\!1}({\cal P})
\]
We thus obtain two compositions 
\[
\begin{array}{c}
\lambda({\cal P}):=(\lambda_1({\cal P}),\lambda_2({\cal P}),\dots)
\\
\lambda'({\cal P}):=(\lambda'_1({\cal P}),\lambda'_2({\cal P}),\dots)
\end{array}
\]
\end{defn}

The following theorem is due to Greene and Fomin, which justifies 
the use of the notation $\lambda'({\cal P})$. 

\begin{thm}
Let ${\cal P}$ be a poset. 
Then $\lambda({\cal P})$ and 
$\lambda'({\cal P})$ are 
partitions. Further, $\lambda'({\cal P})$ is the transpose 
of $\lambda({\cal P})$. 
\end{thm}

Let $x=x_1\cdots x_n$ be a word in $[1,r]^n$, where 
$[1,r]:=\{1,2,\dots,r\}$ is the set of alphabets. 
For a pair $(T,k)$ of a semistandard tableau and 
$k\in[1,r]$, we have the (row) insertion algorithm 
which produces another semistandard tableau. We 
denote this semistandard tableau by $T\leftarrow k$. 

\begin{defn}
Let $x=x_1\cdots x_n\in[1,r]^n$ be a word. 
The semistandard tableau $P(x)$ is defined by
\[
P(x)=\emptyset\leftarrow x_1\leftarrow x_2\leftarrow\cdots
\leftarrow x_n
\]
and is called the $P$-symbol of $x$. 
\end{defn}

\begin{defn}
Let $S_n$ be the symmetric group of $n$ letters acting on 
$[1,n]$. For $w\in S_n$, its permutation poset 
$({\cal P}(w),\le)$ is defined by
\[
\begin{array}{c}
{\cal P}(w)=\{\,(i,w(i))\,|\,i\in[1,n]\,\}\\
\\
(i,w(i))\le(j,w(j))\Leftrightarrow i\le j,\,w(i)\le w(j)
\end{array}
\]
We identify $w\in S_n$ with the word 
$w(1)w(2)\cdots w(n)\in[1,n]^n$. 
\end{defn}

Let $x\in[1,r]^n$ be a word, and assume that $k$ 
appears $n_k$ times in $x$. Then we see $x$ as a 
distinguished coset representative of 
$S_n/S_{n_1}\times\cdots\times S_{n_r}$. 
Thus we can consider permutation posets for arbitrary 
$x\in[1,r]^n$, which we denote by ${\cal P}(x)$. 

\begin{examp}
Let $x=312143\in[1,4]^6$. Then to see it as 
a distinguished coset representative (an element of 
$S_6$) is the same as seeing it as $3_11_121_243_2$. 
Here, we use $1_1\!<\!1_2\!<\!2\!<\!3_1\!<\!3_2\!\!<4$ instead of 
$1\!<\!2\!<\!3\!<\!4\!<\!5\!<\!6$. The permutation 
poset ${\cal P}(x)$ is as follows. 

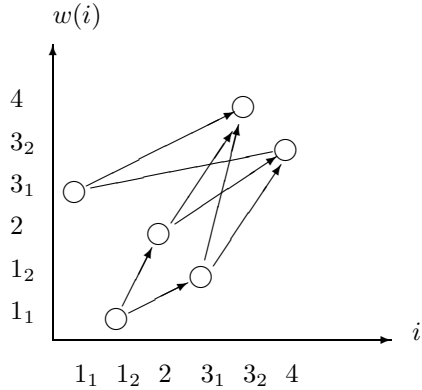
\setlength{\unitlength}{16pt}
\begin{picture}(15,10.5)(-7,-2)
\put(0,0){\vector(1,0){8}}
\put(8.5,0){$i$}
\put(0,0){\vector(0,1){7}}
\put(0,7.5){$w(i)$}
\put(0.5,3.5){\circle{0.5}}
\put(1.5,0.5){\circle{0.5}}
\put(2.5,2.5){\circle{0.5}}
\put(3.5,1.5){\circle{0.5}}
\put(4.5,5.5){\circle{0.5}}
\put(5.5,4.5){\circle{0.5}}
\put(0.5,-1){$1_1$}
\put(1.5,-1){$1_2$}
\put(2.5,-1){$2$}
\put(3.5,-1){$3_1$}
\put(4.5,-1){$3_2$}
\put(5.5,-1){$4$}
\put(-1,0.5){$1_1$}
\put(-1,1.5){$1_2$}
\put(-1,2.5){$2$}
\put(-1,3.5){$3_1$}
\put(-1,4.5){$3_2$}
\put(-1,5.5){$4$}
\put(1.65,0.8){\vector(1,2){0.7}}
\put(1.8,0.65){\vector(2,1){1.4}}
\put(3.8,1.76){\vector(2,3){1.6}}
\put(3.6,1.9){\vector(1,4){0.8}}
\put(2.9,2.75){\vector(3,2){2.4}}
\put(2.8,2.76){\vector(2,3){1.45}}
\put(0.8,3.65){\vector(2,1){3.5}}
\put(0.9,3.6){\line(5,1){4.3}}
\end{picture}

We have $I_1=3,\,I_2=5,\,I_3=6,\,I_4=6,\dots,$ 
and $\lambda({\cal P}(x))=(3,2,1)$. 
\end{examp}

For permutation posets, the following is well known. 

\begin{thm}
\label{RSK shape}
Let $x\in[1,r]^n$ be a word, and ${\cal P}(x)$ be 
its permutation poset. Then 
$\lambda({\cal P}(x))$ equals the shape of the 
$P$-symbol $P(x)$.
\end{thm}

\section{Box-Ball system}

We now recall the box-ball system. Each state is given by 
an infinite sequence of $\{0,1\}$ which has finitely many 
$1$'s. We denote by $1_1,\dots,1_N$ these $1$ read from 
left to right. The description of one step time evolution 
is very simple: for $k=1,\dots,N$, we move $1_k$ to 
the leftmost $0$ among those which sit on the 
right hand side of $1_k$. We give an example. 

\begin{examp}
\[
\begin{array}{cl}
t\,:& \cdots 0\,0\,1\,0\,0\,1\,1\,0\,1\,1\,0\,0\,0\,0\,\cdots\\
\\
t\!+\!1\,:& \cdots 0\,0\,0\,1\,0\,0\,0\,1\,0\,0\,1\,1\,1\,0\,\cdots
\end{array}
\]
It is visualized as follows, and in fact this is the 
original description of the rule. 

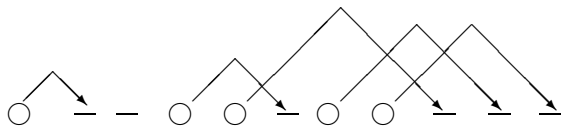
\setlength{\unitlength}{16pt}
\begin{picture}(15,6)(-4,-2)
\put(0.2,0){\circle{0.5}}
 \put(0.3,0.3){\line(1,1){0.7}}
 \put(1.0,1.0){\vector(1,-1){0.8}}
 \put(1.5,0){\line(1,0){0.5}}
 \put(2.5,0){\line(1,0){0.5}}
\put(4.0,0){\circle{0.5}}
 \put(4.3,0.3){\line(1,1){1}}
 \put(5.3,1.3){\vector(1,-1){1.2}}
\put(5.3,0){\circle{0.5}}
 \put(5.6,0.3){\line(1,1){2.2}}
 \put(7.8,2.5){\vector(1,-1){2.4}}
 \put(6.3,0){\line(1,0){0.5}}
\put(7.5,0){\circle{0.5}}
 \put(7.8,0.3){\line(1,1){1.8}}
 \put(9.6,2.1){\vector(1,-1){2}}
\put(8.8,0){\circle{0.5}}
 \put(9.1,0.3){\line(1,1){1.8}}
 \put(10.9,2.1){\vector(1,-1){2}}
\put(10,0){\line(1,0){0.5}}
\put(11.3,0){\line(1,0){0.5}}
\put(12.5,0){\line(1,0){0.5}}
\end{picture}

\end{examp}

For a state, we shall define a finite 
sequence consisting of "$0$", "$($" and "$)$". We first choose 
subsequences of the form "$1\,0$". These are called 
{\bf pairs of stack depth $1$}. We change these pairs 
"$1\,0$" to "$(\,)$". We then delete these pairs from the original 
state, and choose subsequences "$1\,0$" again. 
These pairs are called {\bf pairs of stack depth $2$}. 
We change these pairs "$1\,0$" to "$(\,)$", and delete these 
pairs again. We continue this procedure repeatedly until 
all $1$ are deleted. In the end, all $1$ are made 
into pairs with ")" and their stack depths are defined. 
We say that an opening parenthesis and a closing parenthesis 
are {\bf matched} if they belong to a same pair. 
We now define the stack permutation of the state. 

\begin{defn}
For a state, we associate a finite sequence of "$0$", "(" and ")" 
as above. We read the opening parentheses from left to right, and 
number them $1,2,\dots$ accordingly. 

We number a closing parenthesis $k$ if it makes a pair 
with the $k$ th opening parenthesis. 
The permutation obtained by reading the numbering of 
closing parentheses from left to right is 
called the stack permutation of the state. 
\end{defn}

\begin{examp}
For a state given by
\[
\cdots 0\,0\,1\,0\,0\,1\,1\,0\,1\,1\,0\,0\,0\,0\,\cdots,
\]
we make pairs as follows. 
\[
\begin{array}{c}
\cdots 0\,0\,(\;)\;0\,1\,(\;)\;1\,(\;)\;0\,0\,0\,\cdots\\
\cdots 0\,0\,(\;)\;0\,1\,(\;)\;(\;(\;)\;)\;0\,0\,\cdots\\
\cdots 0\,0\,(\;)\;0\,(\;(\;)\;(\;(\;)\;)\;)\;0\,\cdots
\end{array}
\]
Thus the numbering of the closing parentheses is given by
\[
\begin{array}{c}
\cdots 0\,0\,(\;)\;0\,(\;(\;)\;(\;(\;)\;)\;)\,0\,\cdots\\
\cdots \hphantom{0}\,\hphantom{0}\,\hphantom{0}\;1\,\hphantom{0}\,
\hphantom{0}\;\hphantom{0}\;3\,\hphantom{0}\;
\hphantom{0}\;5\,4\,2\,\hphantom{0}\,\cdots
\end{array}
\]
That is, the stack permutation is $1\,3\,5\,4\,2$. 
\end{examp}

Note that an opening parenthesis does not move to 
the corresponding closing parenthesis, but the total 
set of the opening parentheses moves to the total set of 
the closing parentheses as a whole in each step of 
time evolution. 

We now describe the box-ball 
system as a discrete dynamical system of a path 
on ${\mathbb Z}\times{\mathbb N}$. 
Recall that we have defined 
a finite sequence of "$0$", "(" and ")" for each state. 
We read the sequence from the first "(", and 
associate a walk on ${\mathbb Z}\times{\mathbb N}$ 
starting from $(0,0)$ by the 
rule that "$($" corresponds to "$\uparrow$", and "$0$", "$)$" 
correspond to "$\rightarrow$". 
We see it as a walk from $(-\infty,0)$ to $(\infty,N)$, 
where $N$ is the number of balls in the box-ball system, 
by adding infinitely many "$\rightarrow$" to both sides. 
We give an example. 

\begin{examp}
For a state in the previous example, we have obtained 
the following sequence. 
\[
(\;)\;0\,(\;(\;)\;(\;(\;)\;)\;)
\]
Thus we have a walk as 
follows. The permutation poset 
of the stack permutation $w=13542$ is given by circles. 

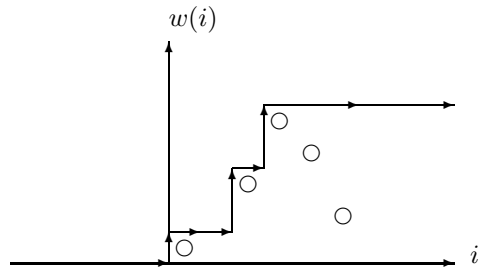
\setlength{\unitlength}{12pt}
\begin{picture}(15,10)(-8,-1)
\put(-1,0){\vector(1,0){9}}\put(8.5,0){$i$}
\put(-1,0){\vector(0,1){7}}\put(-1,7.5){$w(i)$}

\put(-6,0){\vector(1,0){5}}
\put(-1,0){\vector(0,1){1}}
\put(-1,1){\vector(1,0){1}}
\put(0,1){\vector(1,0){1}}
\put(1,1){\vector(0,1){2}}
\put(1,3){\vector(1,0){1}}
\put(2,3){\vector(0,1){2}}
\put(2,5){\vector(1,0){3}}
\put(5,5){\vector(1,0){3}}
\put(-0.5,0.5){\circle{0.5}}
\put(1.5,2.5){\circle{0.5}}
\put(2.5,4.5){\circle{0.5}}
\put(3.5,3.5){\circle{0.5}}
\put(4.5,1.5){\circle{0.5}}

\end{picture}

\end{examp}

Let us think of 
a vertical edge corresponding to "(" of a subsequence 
"0\,(\;". Since the left neighbor of the vertical 
edge is "$0$", we have that 
a positive number of "(" remain waiting to be matched 
as long as the number of "(" is greater than the number of ")". 
It implies that if the vertical edge 
starts at $(x_0,y_0)$, and the walk is in the area 
$y\!-\!y_0>x\!-\!x_0$, then all edges are "(" or ")", and 
no "$0$" appears. 
Using this, we can give a simple evolution rule of the path: 
Let us start walking with $(x_0,y_0)$ and 
continue walking on the path until it hits the line 
$y\!-\!y_0=x\!-\!x_0$ again. We say that 
these edges constitute a group. We partiton the edges 
of the walk into such groups. 
For each group with $(x_0,y_0)$ as above, we reflect the edges in this group 
with respect to the line $y\!-\!y_0=x\!-\!x_0$. 
We then move the zero of the $(x,y)$-plane to the 
lower vertex of the first vertical edge. 
These complete one step time evolution. 
This evolution rule is best understood by an example. 

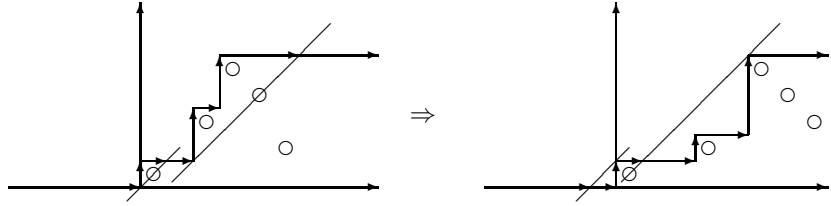
\setlength{\unitlength}{10pt}
\begin{picture}(15,10)(-5,-1)
\put(-1,0){\vector(1,0){9}}
\put(-1,0){\vector(0,1){7}}

\put(-6,0){\vector(1,0){5}}
\put(-1,0){\vector(0,1){1}}
\put(-1,1){\vector(1,0){1}}
\put(0,1){\vector(1,0){1}}
\put(1,1){\vector(0,1){2}}
\put(1,3){\vector(1,0){1}}
\put(2,3){\vector(0,1){2}}
\put(2,5){\vector(1,0){3}}
\put(5,5){\vector(1,0){3}}
\put(-0.5,0.5){\circle{0.5}}
\put(1.5,2.5){\circle{0.5}}
\put(2.5,4.5){\circle{0.5}}
\put(3.5,3.5){\circle{0.5}}
\put(4.5,1.5){\circle{0.5}}
\put(-1.5,-0.5){\line(1,1){2}}
\put(0.2,0.2){\line(1,1){6}}
\put(9.2,2.5){$\Rightarrow$}
\put(15.5,-0.5){\line(1,1){2}}
\put(17.2,0.2){\line(1,1){6}}
\put(17,0){\vector(1,0){8}}
\put(17,0){\vector(0,1){7}}

\put(12,0){\vector(1,0){4}}
\put(16,0){\vector(1,0){1}}
\put(17,0){\vector(0,1){1}}
\put(17,1){\vector(1,0){1}}
\put(18,1){\vector(1,0){2}}
\put(20,1){\vector(0,1){1}}
\put(20,2){\vector(1,0){2}}
\put(22,2){\vector(0,1){3}}
\put(22,5){\vector(1,0){3}}
\put(17.5,0.5){\circle{0.5}}
\put(24.5,2.5){\circle{0.5}}
\put(22.5,4.5){\circle{0.5}}
\put(23.5,3.5){\circle{0.5}}
\put(20.5,1.5){\circle{0.5}}
\end{picture}

Since the closing parentheses of the pairs of "(" and ")" 
define the stack permutation, we have natural 
correspondence between the vertices of the 
permutation poset of the stack permutation and 
the matching pairs: for each vertical edge, 
we choose the horizontal edge in the matching pair. 
If the vertical edge is on the $j$ th row and the 
horizontal edge is on the $i$ th column, then $(i,j)$ is 
an element of the permutation poset. 
In the above example, 
the leftmost three vertices 
correspond to the pairs of stack depth $1$, and 
the left vertex of the remaining two corresponds to 
the pair of stack depth $2$, and the rightmost vertex 
corresponds to the pair of stack depth $3$. 

To describe the poset structure, it is convenient to 
describe the matching pairs by framed boxes. 
For the left example of the above, we have 
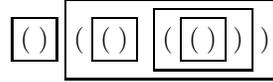
\begin{center}
\fbox{(\;)}\;\fbox{(\;\fbox{(\;)}\;
\fbox{(\;\fbox{(\;)}\;)}\;)}\;
\end{center}

We say that two pairs are in {\bf outer relation} 
if each of the corresponding 
framed box does not contain the other. If one framed box 
contains the other, we say that these are in 
{\bf inner relation}. Note that two pairs are either in 
outer relation or inner relation by the definition of 
the pairs. 

\begin{lem}
\label{lemma 1}
Assume that two pairs of "$($" and "$)$" are in outer 
(resp. inner) relation. Then the corresponding vertices in the 
permutation poset are comparable (resp.not comparable). 
\end{lem}
(Proof) If they are in outer relation, their positions in 
the sequence of "(" and ")" are given by 
\[
\cdots (\;\cdots\; )\;\cdots\; (\;\cdots\;)\;\cdots.
\]
If the left "(" is the $i$ th opening parenthesis and 
the right "(" is the $j$ th opening parenthsis, we have 
$i<j$ and $w(i)<w(j)$ by the definition of the stack 
permutation. Hence the corresponding vertices in 
the permutation poset are comparable. 
The argument for the inner case is similar. 
\QED

\begin{prop}
\label{proposition 1}
To each state, we associate the permutation poset 
of the stack permutation of the state as above. Then 
its vertices corresponding to pairs of stack depth $k$ 
form a chain in the poset. We call it the depth $k$ 
chain and denote it by $C_k$. We then have that 
$|C_1\sqcup\cdots\sqcup C_k|$ gives the maximal number 
of vertices covered by $k$ chains. 
\end{prop}
(Proof) By Lemma \ref{lemma 1}, $C_k$ is a chain. 
We show that this permutation poset admits decomposition 
into disjoint union of antichains $A_k$ such that each $A_k$ 
has the form $\{v_1,\dots,v_{l_k}\}$ where $v_i$ is 
a vertex corresponding to a pair of stack depth $i$. 
Assume that we have already distributed 
vertices of stack depth smaller than $k$ into such 
antichains. In the definition of the pairs, 
it corresponds to the stage that we have deleted 
"$1\,0$"'s $k\!-\!1$ times. By the definition 
of the pairs, each framed box of stack depth $k$ contains a 
framed box of stack depth $k\!-\!1$, and these framed boxes 
of stack depth $k$ are in outer relation. The latter implies 
that we can choose distinct framed boxes of depth $k\!-\!1$ for 
framed boxes of depth $k$. 
Since the vertices of stack depth $k\!-\!1$ are distributed 
to distinct antichains, we can distribute 
the vertices of stack depth $k$ to antichains 
without violating the required property. 

We now assume that $C_1'\cup\cdots\cup C_k'$ gives the 
maximal number of vertices covered by $k$ chains. 
Since each antichain intersects $C_1'\cup\cdots\cup C_k'$ 
at most $k$ times, we can move these vertices into 
$C_1\cup\cdots\cup C_k$ keeping them mutually distict. 
This is possible by the existence of the antichain decomposition 
we have just proved. Hence, $|C_1'\cup\cdots\cup C_k'|$ can not 
exceed $|C_1\cup\cdots\cup C_k|$. 
\QED

\medskip
If we denote by ${\cal P}^t$ the permutation poset at time $t$, 
and by $C_k^t$ the depth $k$ chain of ${\cal P}^t$, 
we have $\lambda_k({\cal P}^t)=|C_k^t|$ by Proposition 
\ref{proposition 1}. 
By Theorem \ref{RSK shape}, we have that $\lambda({\cal P}^t)$ 
is nothing but the shape of the $P$-symbol of the stack 
permutation of the state at time $t$. Hence, the following 
theorem is almost obvious. It simply says that 
the length of depth $k$ chain is conserved, which 
is easily seen from the evolution rule of the path as 
follows. 

\begin{thm}[\cite{TTS}]
For each state, we compute its stack permutation. 
Then the shape of its $P$-symbol is conserved under 
time evolution. 
\end{thm}
(Proof) We show that $|C_k^t|$ is conserved. 
For $k=1$, the elements 
of the chain correspond to convex corners of the path. 
Hence it is obviously conserved by the evolution rule 
of the path. We then delete the depth $1$ chain from 
the posets. This is the same as deleting convex corners 
from the original path and the reflected path. 

To know that the deletion of convex corners 
from the original path and the reflected path 
gives a same walk, it is enough to see that 
deleting convex corners gives the same walk as 
deleting concave corners. To compare the location of $1$'s in 
the walks, We divide the cases by looking at vertical lines 
(the middle lines of the figures below). For the location 
of $0$'s, we divide the cases by looking at horizontal 
lines and the argument is entirely similar, which we omit. 

The leftmost figure represents the case that 
we have vertical lines on both sides. 
One may subdivide the case into four by 
separating the case that 
there is exactly one $0$ in the middle of $1$'s from the case that 
there are more than one $0$'s in the 
middle of $1$'s, if one wishs. The remaining two cases are the left end and 
the right end of the walk. 

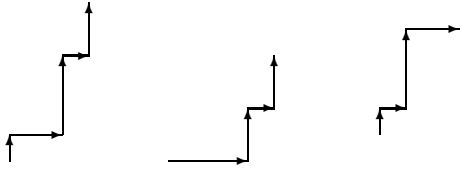
\setlength{\unitlength}{10pt}
\begin{picture}(15,9)(-7,-1)

\put(0,0){\vector(0,1){1}}
\put(0,1){\vector(1,0){2}}
\put(2,1){\vector(0,1){3}}
\put(2,4){\vector(1,0){1}}
\put(3,4){\vector(0,1){2}}

\put(6,0){\vector(1,0){3}}
\put(9,0){\vector(0,1){2}}
\put(9,2){\vector(1,0){1}}
\put(10,2){\vector(0,1){2}}

\put(14,1){\vector(0,1){1}}
\put(14,2){\vector(1,0){1}}
\put(15,2){\vector(0,1){3}}
\put(15,5){\vector(1,0){2}}
\end{picture}

By comparing 
the results of the deletion of the concave corners and 
the convex corners, we know that the new walks are the same. 
In particular, deleting convex corners from 
the original path and the 
reflected path give a same walk. 
Since the depth $2$ chain becomes 
the depth $1$ chain of the new poset, we can apply 
the same argument to conclude that 
$|C_2^t|$ is conserved. By repeated use of the argument, 
we also have the conservation laws for all $k$. \QED

\section{Energy functions and stack permutations}

We now turn to the crystal description of the box-ball 
system. Let $B:=\{\fbox{0},\,\fbox{1}\,\}$ 
be the $A_1$ crystal associated 
with the vector representation 
whose highest weight vector is $\fbox{0}$. Its affinization 
is denoted by $A\!f\!f(B):={\mathbb Z}\times B$. This is an 
$A_1^{(1)}$ crystal. Note that the numbering of $\fbox{0}$ and 
$\fbox{1}$ is different from the usual one. 
$B$ is identified with the subset 
$\{0\}\times B$. For each state, 
we cut sufficiently remote $0$'s and consider it as 
an element in $B^{\otimes n}$. 

To describe the time evolution rule, we take the crystal of 
the $l$ th symmetric tensor $B_l$ and its affinization 
$A\!f\!f(B_l)$ with $l$ sufficiently large. The elements 
of $B_l$ are nondecreasing sequences of length $l$ 
whose entries are $0$ and $1$. We write $0^{m_1}1^{m_2}$ 
$(m_1\!+\!m_2=l)$ for these elements. 
We use combinatorial $R$ matrices to get isomorphism 
of affine crystals as follows. 

\setlength{\unitlength}{8pt}
\begin{picture}(20,10)(-6,-5)
\put(0,0.3){\vector(1,0){23}}
\put(-5,0){$A\!f\!f(B_l)$}
\put(23.5,0){$A\!f\!f(B_l)$}
\put(-1.5,2){$\otimes$}
\put(22.3,-2){$\otimes$}

\put(1,3){$Af\!f(B)$}
\put(3,2){\vector(0,-1){4}}
\put(1,-3.5){$A\!f\!f(B)$}

\put(5.5,3){$\otimes$}
\put(5.5,-3.5){$\otimes$}

\put(7,3){$Af\!f(B)$}
\put(9,2){\vector(0,-1){4}}
\put(7,-3.5){$A\!f\!f(B)$}

\put(11.5,3){$\otimes$}
\put(11.5,-3.5){$\otimes$}
\put(13.5,3){$\cdots$}
\put(13.5,-3.5){$\cdots$}
\put(15.5,3){$\otimes$}
\put(15.5,-3.5){$\otimes$}
\put(17,3){$Af\!f(B)$}
\put(19,2){\vector(0,-1){4}}
\put(17,-3.5){$A\!f\!f(B)$}
\end{picture}

After we embed $B^{\otimes n}$ to $A\!f\!f(B)^{\otimes n}$, 
we apply this combinatorial 
row-to-row transfer to the tensor product of 
\fbox{$0^l$} with the upper $A\!f\!f(B)^{\otimes n}$ 
to get the lower $A\!f\!f(B)^{\otimes n}$ tensored by 
\fbox{$0^l$}. Then we forget the symmetric tensor part and 
the ${\mathbb Z}$ part of the affine crystal. 
The result is an element 
of $B^{\otimes n}$. This procedure gives one step 
time evolution of the box-ball system. 

We consider the isomorphism for arbitrary $l$. Then 
for a state $p$, we have 
\[
\left(0\times \fbox{$0^l$}\right)\otimes p \mapsto 
p' \otimes\left(E_l(p)\times \fbox{$0^{m_1}1^{m_2}$}\right)
\]
for some $p'\in A\!f\!f(B)^{\otimes n}$ and 
$m_1,m_2\in{\mathbb N}$. 
These $E_l(p)$ are called energy functions. It is known 
\cite{NY} that if we set $E_l=0$ and increase it 
by one at the sites of the following form, then 
the final value of $E_l$ coincides with $E_l(p)$. 

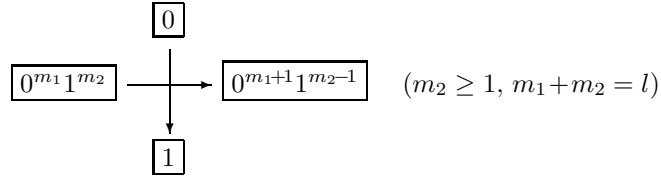
\setlength{\unitlength}{8pt}
\begin{picture}(20,10)(-12,-5)
\put(0,0.3){\vector(1,0){4}}
\put(-5.5,0){\fbox{$0^{m_1}1^{m_2}$}}
\put(4.5,0){\fbox{$0^{m_1\!+\!1}1^{m_2\!-\!1}$}}
\put(1.2,3){\fbox{$0$}}
\put(2,2){\vector(0,-1){4}}
\put(1.2,-3.5){\fbox{$1$}}
\put(13,0){$(m_2\ge1,\,m_1\!+\!m_2=l)$}
\end{picture}

By using the fact that time evolution is obtained from 
a crystal isomorphism 
of affine crystals, Fukuda, Okado and Yamada 
\cite[Theorem 3.2]{FOY} have proved that 
these $E_l(p)$ are conserved quantities of this box-ball system. 

The purpose of this section is to relate these quantities to 
the stack permutation of the state $p$. 

\begin{thm}
For each state $p$, we define the sequence of "$0$", "(" and ")" as 
in the previous section. Then $E_l$ increases 
precisely at the sites corresponding to ")". 
\end{thm}
(Proof) Assume that ")" corresponds to a pair of depth 
$k$. We shall show that if $l\ge k$, $E_l(p)$ does 
increase at this site. 

Let "(" be the corresponding opening parenthesis, and 
$1^{k_1}0^{l_1}\cdots 1^{k_N}0^{l_N}$ be the walk starting 
from the vertical edge corresponding to the "(" and ending at 
the horizontal edge corresponding to the ")". We write the 
evolution of the symmetric tensor along the path as follows. 

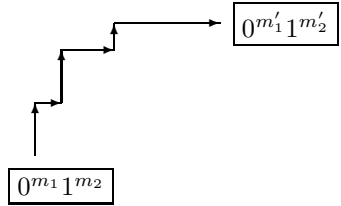
\setlength{\unitlength}{10pt}
\begin{picture}(20,10)(-12,-2)
\put(-1,-0.5){\fbox{$0^{m_1}1^{m_2}$}}
\put(0,1){\vector(0,1){2}}
\put(0,3){\vector(1,0){1}}
\put(1,3){\vector(0,1){2}}
\put(1,5){\vector(1,0){2}}
\put(3,5){\vector(0,1){1}}
\put(3,6){\vector(1,0){4}}
\put(7.5,5.5){\fbox{$0^{m_1'}1^{m_2'}$}}

\end{picture}

Assume that 
$E_l$ does not increase at this site. Then the status 
of the symmetric tensor on 
the both ends of the last edge is \fbox{$0^l$}. 
We denote by $e_1,\dots,e_N$ the last edges of 
$1^{k_1},\dots,1^{k_N}$ respectively. 
Then we can prove the following by downward induction 
on $i$. 

\begin{itemize}
\item
The symmetric tensor 
on the upper end of $e_i$ has the form 
\fbox{$0^{l-s}1^s$} with $s<l$. 
\item
Since the upper end of $e_i$ is not saturated, 
we have steady increase of the number of 
$1$ in the symmetric tensor during $k_i$ 
vertical edges, and no saturation occurs during these edges. 
\item
Since $l_{i\!-\!1}\!-\!k_i\!+\!\cdots\!-\!k_N\!+\!l_N\le k$, 
the symmetric tensor 
on the upper end of $e_{i\!-\!1}$ has the form 
\fbox{$0^{l-s}1^s$} with $s<l$. 
\end{itemize}

We then have the following. 

\begin{itemize}
\item
Since $k_1\!-\!l_1\!+\!\cdots\!-\!l_i\ge0$, 
we have steady decrease of the number of 
$1$ in the symmetric tensor 
during $l_i$ horizontal edges, 
and no saturation occurs. 
\end{itemize}

Therefore, we conclude that the left end of the last 
edge of the walk has the symmetric tensor of the form 
\fbox{$0^{l-s}1^s$} with $s>\sum k_i\!-\!\sum l_i\ge0$, which 
contradicts the assumption at the beginning. 
(In particular, we have that $I_l({\cal P})\le E_l(p)$ where 
${\cal P}$ is the permutation poset of the stack 
permutation of the state $p$.) 

We nextly show that if $l<k$, then $E_l$ does not 
increase at this site. To prove this, we show for 
arbitrary $k,l$ that 
the right end of the last edge of the walk has the 
symmetric tensor \fbox{$0^l$} if the stack depth $k$ is 
equal or greater 
than $l$, and \fbox{$0^{m_1}1^{m_2}$} $(m_1\ge k)$ if 
$k$ is equal or smaller than $l$. We prove it by 
induction on $k$. 
If $k=1$, the proof is obvious. If 
$k\le l$, 
we choose the last closing parethesis of stack depth 
$k\!-\!1$. Then by the induction hypothesis, the 
symmetric tensor has the form 
\fbox{$0^{m_1}1^{m_2}$} $(m_1\ge k\!-\!1)$ at this site. 
Note that we have already proved that no saturation 
occurs during the walk if $k\le l$. 
Hence, if we start the walk with 
\fbox{$0^{m_1}1^{m_2}$}, we end the walk with \fbox{$0^{m_1}1^{m_2}$}. 
From this, we know that the symmetric tensor at the 
left end of the last edge of the walk is also 
\fbox{$0^{m_1}1^{m_2}$} $(m_1\ge k\!-\!1)$. 
Hence, the right end of the last edge has the form 
\fbox{$0^{m_1'}1^{m_2'}$} $(m_1'\ge k)$ if $m_2>0$ and 
\fbox{$0^l$} if $m_2=0$. But we also have $l\ge k$ in the latter case. 

We now assume that $k\ge l$. We choose the last closing 
parenthesis of stack depth equal or greater than $l$. Since its 
stack depth is smaller than $k$, we can apply the induction 
hypothesis to know that the symmetric tensor has the form 
\fbox{$0^l$} at this site. Further, 
since we have pairs of stack depth less than $l$ 
during this site and the last edge of the walk, 
we have that the symmetric tensor at the 
left end of the last edge has \fbox{$0^l$}. Thus the same is true 
for the right end of the last edge. 

Therefore, we have proved that 
$E_l$ increases precisely at the sites corresponding to 
")". (In particular, we have also proved that $I_l({\cal P})=E_l(p)$.) 
\QED

\section{Conclusion}

For a state $p^t$ at time $t$, we denote by ${\cal P}^t$ the 
permutation poset of the stack permutation of the state 
$p^t$. Then the energy function $E_l(p^t)$ counts the 
vertices of ${\cal P}^t$ whose 
stack depth are equal or less than $l$. On the other 
hand, the $l$ th row of the shape $\lambda(p^t)$ of the $P$-symbol 
of the stack permutation is equal to the number of 
vertices of ${\cal P}^t$ whose stack depth are $l$. 
Hence these quantities are naturally explained by the notion 
of stack depth, and we have $I_l({\cal P}^t)=E_l(p^t)$. 

Further, the evolution rule of a path naturally explains 
why these quantities are conserved.

\bigskip
Tokyo University of Mercantile Marine, 

\smallskip
Tokyo 135-8533, Japan.

\end{document}